\documentclass[]{amsart}

\usepackage{amsmath, amsthm, amssymb}
\usepackage{amsfonts}
\usepackage{graphicx}
\usepackage{wrapfig}

\newtheorem{Thm}{Theorem}[section]
\newtheorem{Lem}[Thm]{Lemma}

\theoremstyle{definition}
\newtheorem{Def}[Thm]{Definition}

\begin{document}    
\title{A genus zero Lefschetz fibration on the Akbulut cork} 
\author{Takuya Ukida}
\address{Department of Mathematics, Tokyo Institute of Technology,
2-12-1 Oh-okayama, Meguro-ku, Tokyo 152-8551, Japan}
\email{ukida.t.aa@m.titech.ac.jp}
\date{}
\thanks{The author was partially supported by Research Fellow of Japan
Society for the Promotion of Science}
\begin{abstract}
    We first construct a genus zero positive allowable Lefschetz fibration 
    over the disk (a genus zero PALF for short) on the Akbulut cork 
    and describe the monodromy as a positive factorization 
    in the mapping class group of a surface of genus zero with five boundary components. 
    We then construct genus zero PALFs on infinitely many exotic pairs of 
    compact Stein surfaces such that one is a cork twist of the other along an Akbulut cork. 
    The difference of smooth structures on each of exotic pairs of compact Stein surface 
    is interpreted as the difference of the corresponding 
    positive factorizations in the mapping class group 
    of a common surface of genus zero. 
\end{abstract}
\maketitle

    \section{Introduction}
    Gompf \cite{Go} proved that compact Stein surfaces can be characterized in terms of 
    handle decompositions, or more precisely, Kirby diagrams. 
    Akbulut and Yasui \cite{AY1} introduced corks and plugs, 
    which are compact Stein surfaces themselves, and constructed 
    various exotic smooth structures on Stein surfaces by using cork twists and plug twists 
    together with Gompf's characterization and Seiberg-Witten invariants. On the other hand, 
    Loi and Piergallini \cite{LP} proved that every compact Stein surface admits a positive allowable 
    Lefschetz fibration over $D^2$ (a PALF for short), which enables us to investigate compact 
    Stein surfaces in terms of positive factorizations in mapping class groups 
    (see also Akbulut and Ozbagci \cite{AO}, Akbulut and Arikan \cite{AA}). 

\begin{figure}[htbp]
        \centering
        \includegraphics[width=20mm]{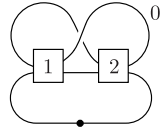}
        \caption{Kirby diagram for $W_1$.}
        \label{AkbulutCork}
    \end{figure}

    In this paper we first construct a genus zero PALF on the Akbulut cork 
    and describe the monodromy as a positive factorization 
    in the mapping class group of a fiber. 
    The Akbulut cork is the pair $(W_1,f_1)$ of the manifold $W_1$ shown
    in Figure \ref{AkbulutCork}
    and an involution $f_1$ on $W_1$ (see Definition \ref{def:AkbulutCork}). 
    The manifold $W_1$ is often called the Mazur manifold.

\begin{figure}[htbp]
        \centering
        \includegraphics[width=70mm]{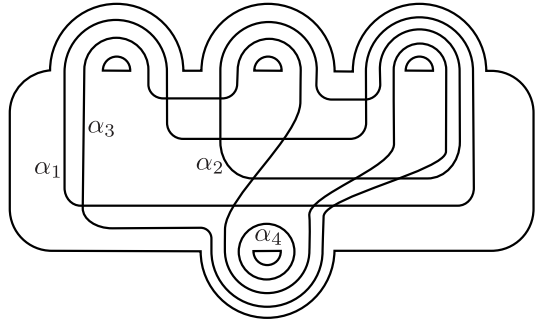}
        \caption{Vanishing cycles of a genus zero PALF on $W_1$.}
        \label{PALF_of_the_Akbulut_cork}
    \end{figure}

    \begin{Thm}\label{main1}
        The manifold $W_{1}$
        admits a genus zero PALF.
        The monodromy of the PALF is described by the
        factorization
        $t_{\alpha_4} t_{\alpha_3} t_{\alpha_2} t_{\alpha_1}$, 
        where $t_\alpha$ is a right-handed Dehn twist along a
        simple closed curve $\alpha$ on a fiber and
        $\alpha_4, \ldots, \alpha_1$ are simple closed curves
        shown in Figure \ref{PALF_of_the_Akbulut_cork}.
    \end{Thm}

    Note that the genus of a PALF on the manifold $W_1$
    which is obtained by applying any known method
    (cf. \cite{AO} and \cite{AA} ) is much larger than zero.

\begin{figure}[htbp]
        \centering
        \includegraphics[width=30mm]{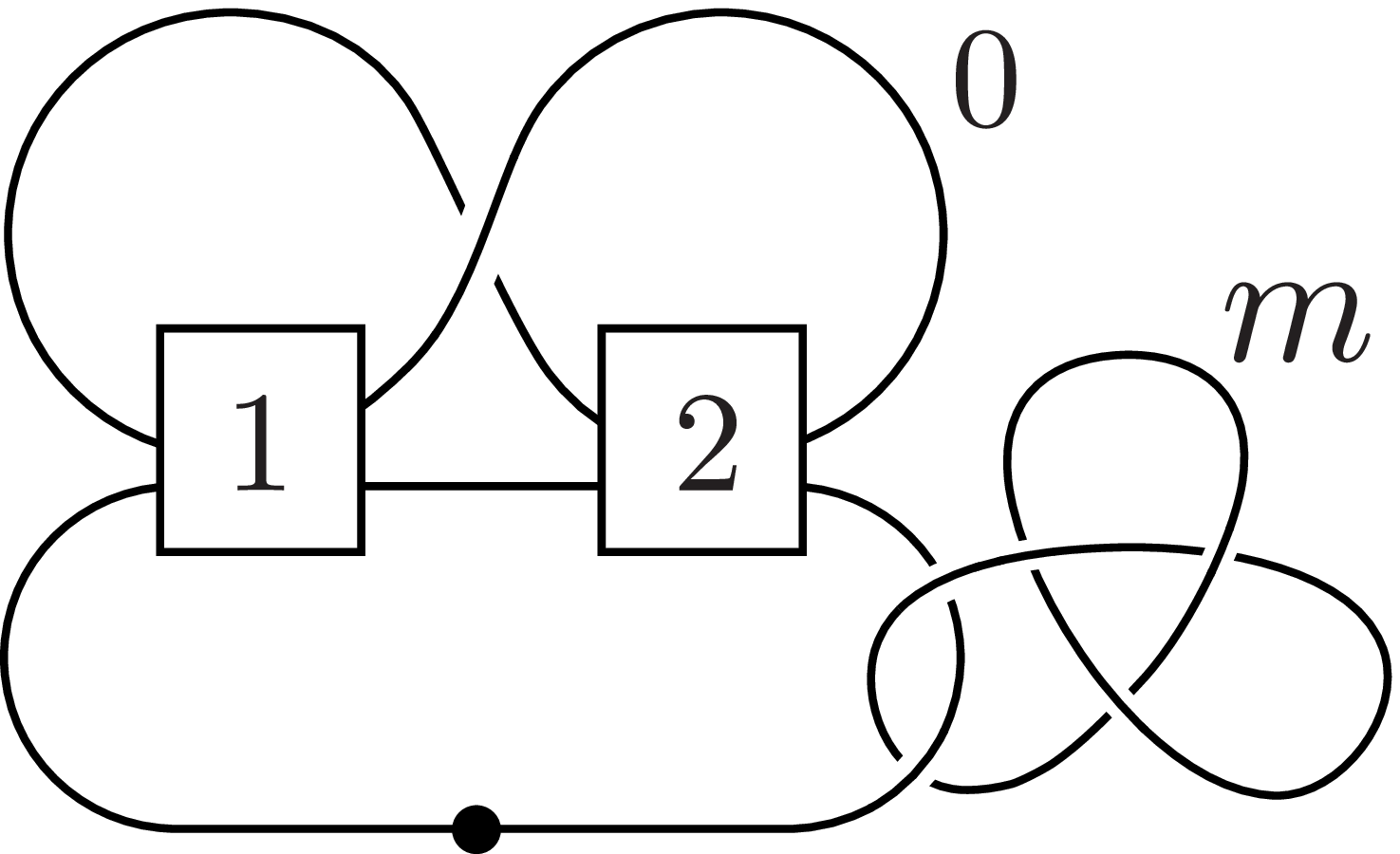}
        \caption{Kirby diagram for $C_1(m, 1, 3, 0)$.}
        \label{Kirby_diag_of_C1}
    \end{figure}

\begin{figure}[htbp]
        \centering
        \includegraphics[width=30mm]{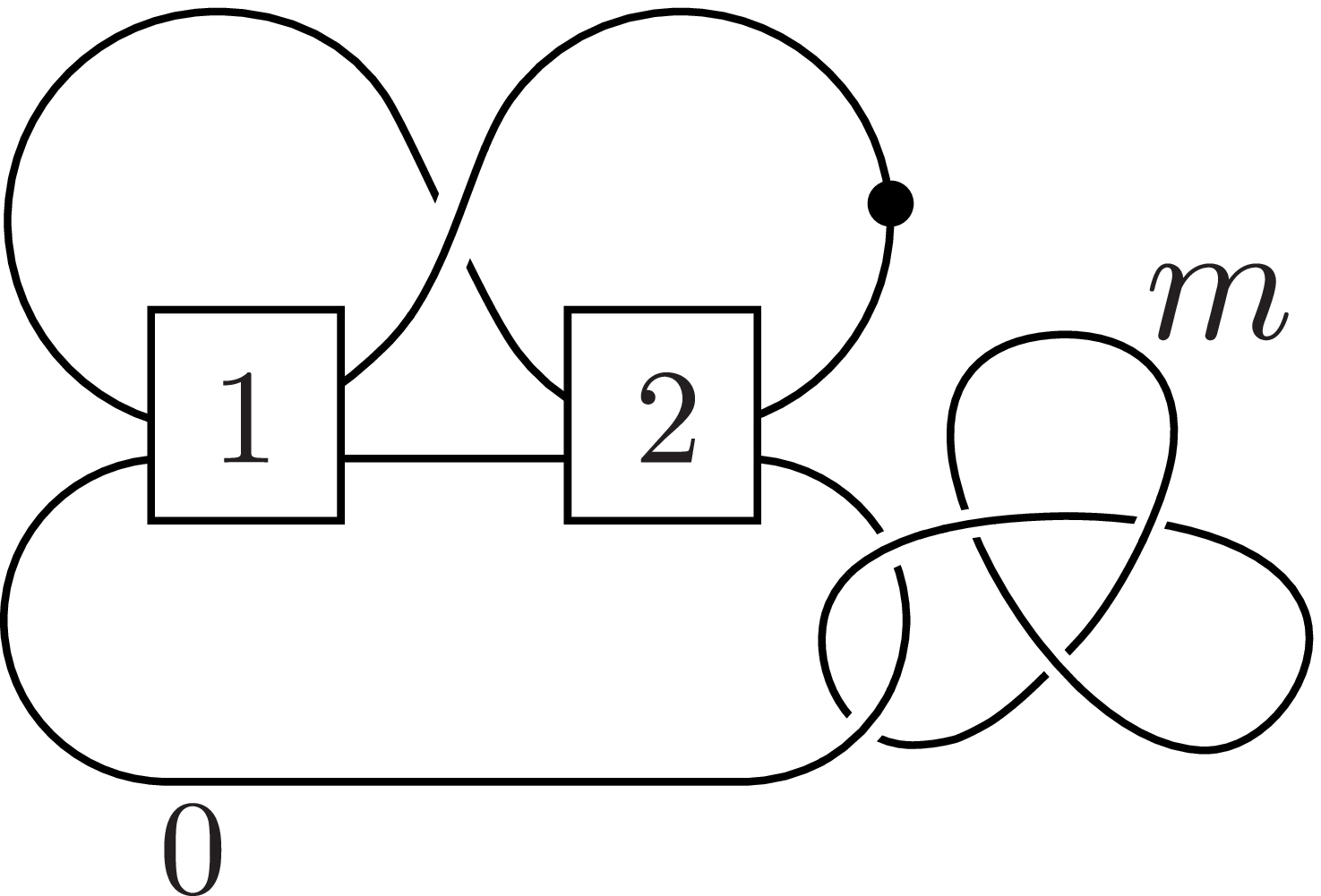}
        \caption{Kirby diagram for $C_2(m, 1, 3, 0)$.}
        \label{Kirby_diag_of_C2}
    \end{figure}

    Akbulut and Yasui \cite{AY2} proved that the compact Stein surfaces 
    $C_1(m,1,3,0)$ and $C_2(m,1,3,0)$ $(m\leq -5)$ 
    shown in Figure \ref{Kirby_diag_of_C1} and Figure \ref{Kirby_diag_of_C2}
    are homeomorphic but not diffeomorphic to each other. 
    It is easily seen that $C_2(m,1,3,0)$ is a cork twist of $C_1(m,1,3,0)$ along an obvious Akbulut cork. 

    We next construct PALFs with the same fiber on 
    $C_1(m,1,3,0)$ and $C_2(m,1,3,0)$ for each integer $m$ less than $-4$. 
    The common fiber is a surface of genus zero with $-m+5$ boundary components.

\begin{figure}[htbp]
        \centering
        \includegraphics[width=110mm]{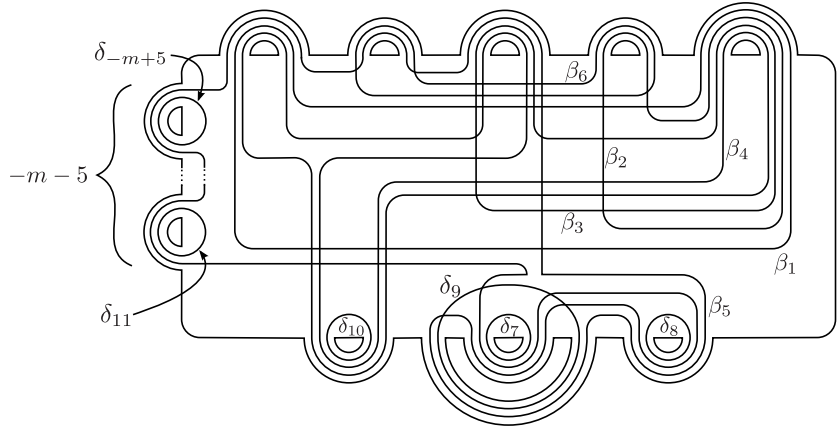}
        \caption{Vanishing cycles of a genus zero PALF on $C_1(m, 1, 3, 0)$.}
        \label{PALF_of_C1}
    \end{figure}

\begin{figure}[htbp]
        \centering
        \includegraphics[width=110mm]{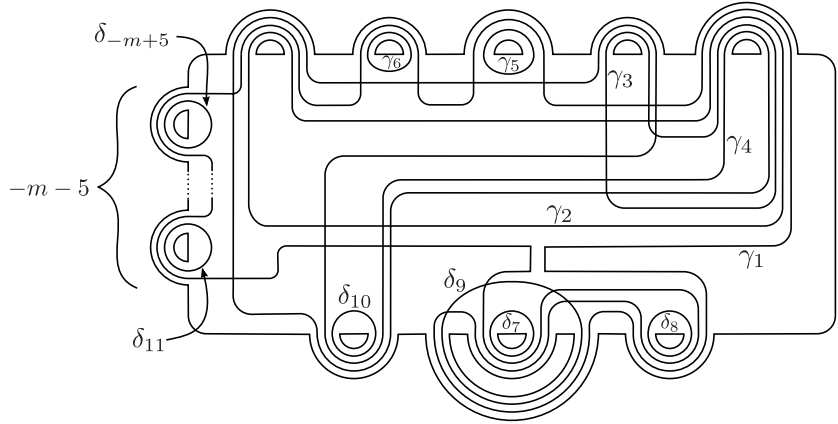}
        \caption{Vanishing cycles of a genus two PALF on $C_2(m, 1, 3, 0)$.}
        \label{PALF_of_C2}
    \end{figure}

    \begin{Thm}\label{main2}
        The manifolds $C_1(m, 1, 3, 0)$ and $C_2(m, 1, 3, 0)$ $(m \leq -5)$
        shown in Figure \ref{Kirby_diag_of_C1} and  Figure \ref{Kirby_diag_of_C2}
        admit genus zero PALFs.
        The monodromy of the PALF on $C_1(m, 1, 3, 0)$ is
        described by the positive factorization
        \begin{align*}
            t_{\delta_{-m+5}} \ldots t_{\delta_{11}} t_{\delta_{10}}
                t_{\delta_9} t_{\delta_8} t_{\delta_7}
                t_{\beta_6} t_{\beta_5} t_{\beta_4} t_{\beta_3}
                t_{\beta_2} t_{\beta_1}
        \end{align*}
        while that for $C_2(m, 3, 1, 0)$ is described by the
        positive factorization
        \begin{align*}
            t_{\delta_{-m+5}} \ldots t_{\delta_{11}} t_{\delta_{10}}
                t_{\delta_9} t_{\delta_8} t_{\delta_7}
                t_{\gamma_6} t_{\gamma_5} t_{\gamma_4} t_{\gamma_3}
                t_{\gamma_2} t_{\gamma_1}
        \end{align*}
        , where $\beta_i, \gamma_j$ are simple closed curves shown in
        Figure \ref{PALF_of_C1} and Figure \ref{PALF_of_C2}.
    \end{Thm}

The difference of smooth structures on $C_1(m,1,3,0)$ and $C_2(m,1,3,0)$ 
(or the effect of cork twisting the former to obtain the latter) 
is reflected in the corresponding 
positive factorizations as the difference between partial factorizations 
$t_{\beta_6}t_{\beta_5}t_{\beta_4}t_{\beta_3}t_{\beta_2}t_{\beta_1}$ and 
$t_{\gamma_6}t_{\gamma_5}t_{\gamma_4}t_{\gamma_3}t_{\gamma_2}t_{\gamma_1}$.

    In Section 2 we briefly review definitions of Mapping class groups,
    PALF, Stein surfaces, and corks, and recall several known results.
    We prove our main theorems in Section 3.

    \section{Preliminaries}

        \subsection{Mapping class groups}

        In this subsection we review a precise definition of
        mapping class groups of surfaces with boundary and that of
        Dehn twists along simple closed curves on surfaces.

        \begin{Def}
            Let $F$ be a compact oriented  connected surface with boundary.
            Let Diff${}^+(F, \partial F)$ be
            the group of all orientation-preserving self-diffeomorphisms
            of $F$ fixing the boundary $\partial F$ point-wise.
            Let Diff${}^+_0(F, \partial F)$ be
            the subgroup of Diff${}^+(F, \partial F)$ consisting of
            self-diffeomorphisms isotopic to the identity.    
            The quotient group
            Diff${}^+(F, \partial F)/$ Diff${}^+_0(F, \partial F)$
            is called the mapping class group of $F$ and it is
            denoted by Map$(F, \partial F)$.
        \end{Def}

        \begin{Def}
            A {\em positive $($or right-handed$)$ Dehn twist}
            along a simple closed curve $\alpha$, $t_\alpha:F \rightarrow F$
            is a diffeomorphism  obtained by
            cutting $F$ along $\alpha$, twisting $360^\circ$
            to the right and regluing.
        \end{Def}

        %%%%%%%%%%%%%%%%%%%%%%%%%%%%%%%%%%%%%%%%%%%%%%%%%%

        \subsection{PALF}

        \begin{Def}\label{aLFdef}
            Let $M^4$ and $B^2$ be compact oriented smooth manifolds
            of dimensions $4$ and $2$.
            Let $f: M \rightarrow B$ be a smooth map.
            $f$ is called a {\em positive Lefschetz fibration} over $B$
            if it satisfies the following conditions (1) and (2):

            \begin{enumerate}
                \renewcommand{\labelenumi}{(\arabic{enumi})}
                \setlength{\itemsep}{0cm} 
                                \setlength{\parskip}{1mm} 
                \item There are finitely many critical values
                      $b_1, \ldots, b_m $ of $f$ in the interior of $B$
                      and there is a unique critical point $p_i$
                      on each fiber $f^{-1}(b_i)$, and
                \item The map $f$ is locally written as
                      $f(z_1, z_2) = z_1^2 + z_2^2$ with respect to
                      some local complex coordinates around
                      $p_i$ and $b_i$ compatible with the orientations of
                      $M$ and $B$.
            \end{enumerate}    
        \end{Def}

        \begin{figure}[htbp]
            \centering
            \includegraphics[width=80mm]{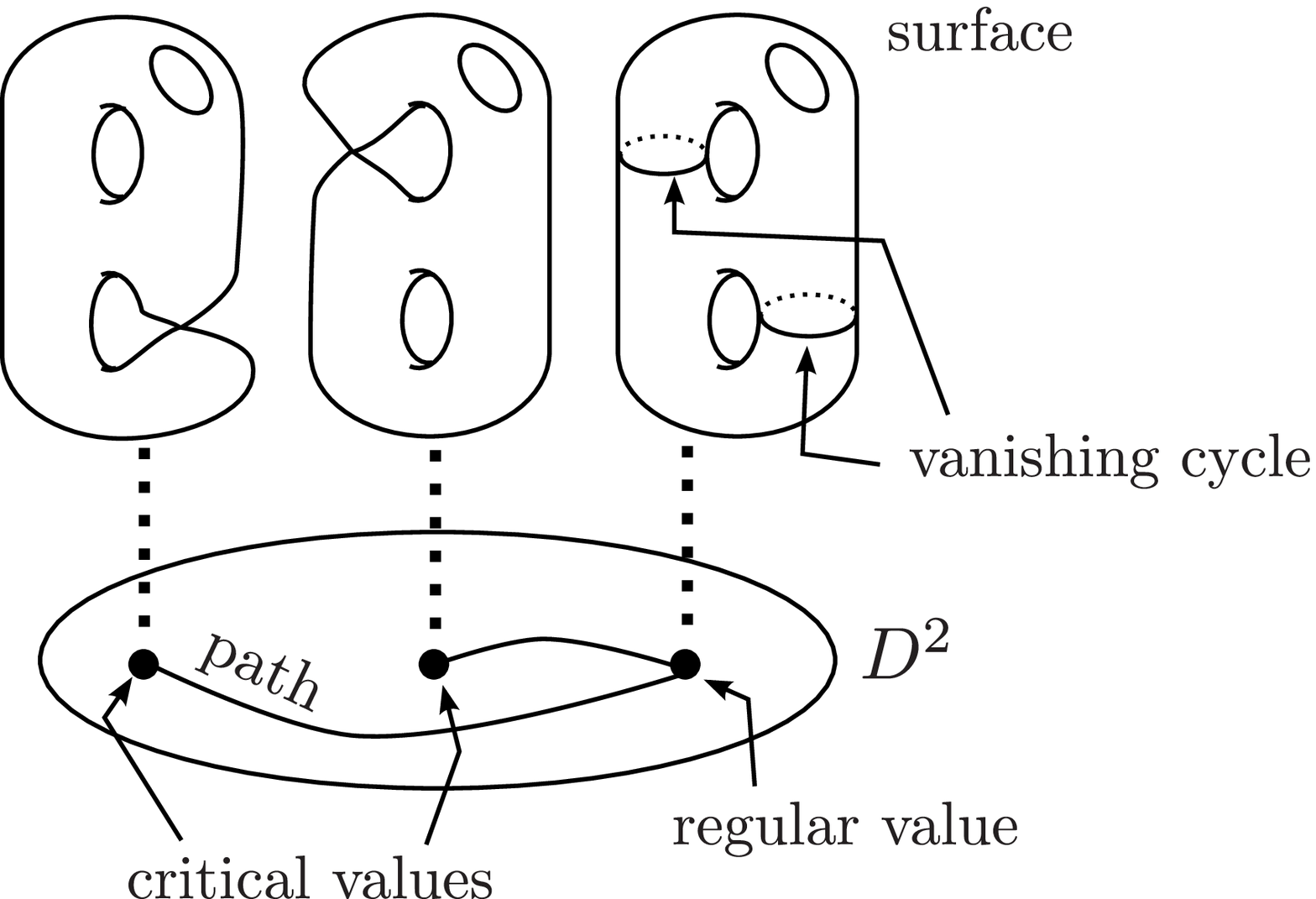}
            \caption{PALF.}
            \label{PALFdiag}
        \end{figure}

        \begin{Def}
            A positive Lefschetz fibration is called {\em allowable} if 
            its all vanishing cycles are homologically non-trivial on the
            fiber.
            A positive allowable Lefschetz fibration over $D^2$ with
            bounded fibers is called a {\em PALF} for short.
        \end{Def}

        The following Lemma is useful to prove Theorem \ref{main1}.

        \begin{Lem}[{cf. Akbulut-Ozbagci \cite[Remark 1]{AO}}]
            \label{Lefschetz_2_h_PALF}
            Suppose that a $4$-manifold $X$ admits a PALF.
            If a $4$-manifold $Y$ is obtained from $X$ by attaching a
            Lefschetz $2$-handle, then $Y$ also admits a PALF.
        \end{Lem}

        The Lefschetz $2$-handle is defined as follows.

        \begin{Def}
            Suppose that $X$ admits a PALF.
            A {\em Lefschetz $2$-handle} is a $2$-handle attached
            along a homologically non-trivial simple closed curve
            in the boundary of $X$ with framing $-1$ relative to
            the product framing induced by the fiber structure. 
        \end{Def}

        %%%%%%%%%%%%%%%%%%%%%%%%%%%%%%%%%%%%%%%%%%%%%%%%%%

        \subsection{Stein surfaces}

        In this section, we recall a definition of Stein surfaces.
        The question of which smooth $4$-manifolds admit
        Stein structures can be completely reduced to a
        problem in handlebody theory.

        \begin{Def}
            A complex manifold is called a {\em Stein manifold} if
            it admits a proper biholomorphic embedding to $\mathbb{C}^n$. 
        \end{Def}

        \begin{Def}
            Let $W$ be a compact manifold with boundary. 
            The manifold $W$ is called a {\em Stein domain} if it
            satisfies following condition:
            There is a Stein manifold $X$ and
            a plurisubharmonic function
            $\varphi : X \rightarrow [0, \infty)$
            such that 
            $W = \varphi^{-1}([0, a])$
            for a regular value $a$ of $\varphi$.
        \end{Def}

        \begin{Def}
            A Stein manifold or a Stein domain is called
            a {\em Stein surface} if its complex dimension is $2$.
        \end{Def}

        %%%%%%%%%%%%%%%%%%%%%%%%%%%%%%%%%%%%%%%%%%%%%%%%%%

        \subsection{Corks}

        Corks are Stein surfaces and
        they are useful for constructing exotic manifolds.

        \begin{Def}
            Let $C$ be a  Stein domain.
            Let $\tau : \partial C \rightarrow \partial C$ be an
            involution on the boundary $\partial C$ of $C$.
            \begin{enumerate}
                \renewcommand{\labelenumi}{(\arabic{enumi})}
                \setlength{\itemsep}{0cm} 
                                \setlength{\parskip}{1mm}
                \item $(C, \tau)$ is called a {\em cork} if $\tau$ 
                      extends to a self-homeomorphism of $C$,
                      but does not extend to any self-diffeomorphism of $C$.
                \item Suppose that $C$ is embedded in a smooth
                      $4$-manifold $X$.
                      The manifold obtained from $X$ by removing $C$
                      and regluing it via $\tau$
                      is called a {\em cork twist} of
                      $X$ along $(C, \tau)$.
                \item The pair $(C, \tau)$ is called a {\em cork of $X$}
                      if the cork twist of $X$ along $(C, \tau)$ is
                      homeomorphic but not diffeomorphic to $X$. 
            \end{enumerate}    
        \end{Def}

        In this paper, we investigate Akbulut cork $(W_1, f_1)$ ([Ak]).

        \begin{Def}\label{def:AkbulutCork}
            Let $W_1$ be a smooth 4-manifold
            given by Figure \ref{AkbulutCork}.
            Let $f_1: \partial W_1 \rightarrow \partial W_1$
            be the obvious involution obtained from first surgering
            $S^1 \times D^3$ to $D^2 \times S^2$ in the
            interiors of $W_1$, then surgering the
            other imbedded $D^2 \times S^2$ back to
            $S^1 \times D^2$.
        \end{Def}

        \begin{Thm}[Akbulut \cite{Ak}]
            The pair $(W_1, f_1)$ is
            a cork.
        \end{Thm}

        Akbulut and Yasui constructed infinitely many exotic pairs.

        \begin{Thm}[{Akbulut-Yasui \cite[Theorem 3.3(2)]{AY2}}]
            $C_1(m, n, p, 0)$ and $C_2(m, n, p, 0)$ are
            homeomorphic but not diffeomorphic to each other,
            for $1 \leq n \leq 3, p \geq 3$ and
            $m \leq p^2 - 3p + 1$.
        \end{Thm}

        \begin{figure}[htbp]
            \centering
            \includegraphics[width=80mm]{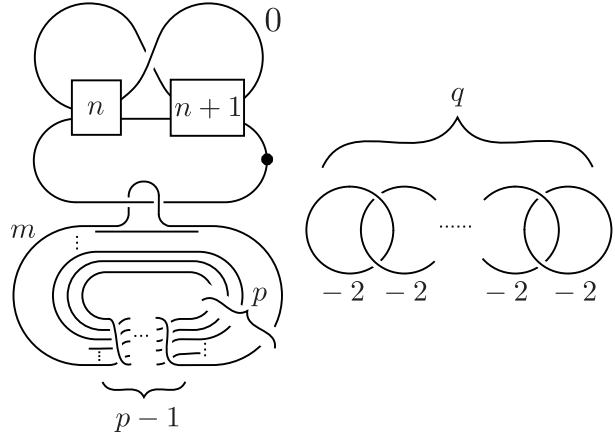}
            \caption{Kirby diagram for $C_1(m, n, p, q)$.}
            \label{C1(mnpq)}
        \end{figure}

        \begin{figure}[htbp]
            \centering
            \includegraphics[width=80mm]{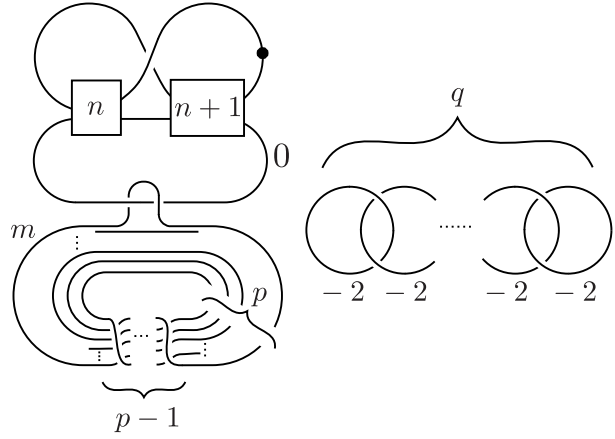}
            \caption{Kirby diagram for $C_2(m, n, p, q)$.}
            \label{C2(mnpq)}
        \end{figure}

        %%%%%%%%%%%%%%%%%%%%%%%%%%%%%%%%%%%%%%%%%%%%%%%%%%

    \section{Proofs of Theorems \ref{main1} and \ref{main2}.}

    {\it Proof of Theorem  \ref{main1}.}
    Let $F$ be the compact oriented surface of genus zero
    with $5$ boundary components
    and $\alpha_1,\ldots , \alpha_4$ be the curves on $F$ shown in
    Figure \ref{PALF_of_the_Akbulut_cork}. 
    We denote the right-handed Dehn twists along $\alpha_1,\ldots , \alpha_4$ 
    by $t_{\alpha_1}, \ldots ,t_{\alpha_4}$, respectively. 
    Let $f:X \rightarrow D^2$ be a Lefschetz fibration over $D^2$ 
    with monodromy representation $(t_{\alpha_4},\ldots ,t_{\alpha_1})$. 
    Since each curve $\alpha_i$ is homologically non-trivial on $F$, 
    we see that $f$ is a PALF with fiber $F$. 

\begin{figure}[htbp]
        \centering
        \includegraphics[width=120mm]{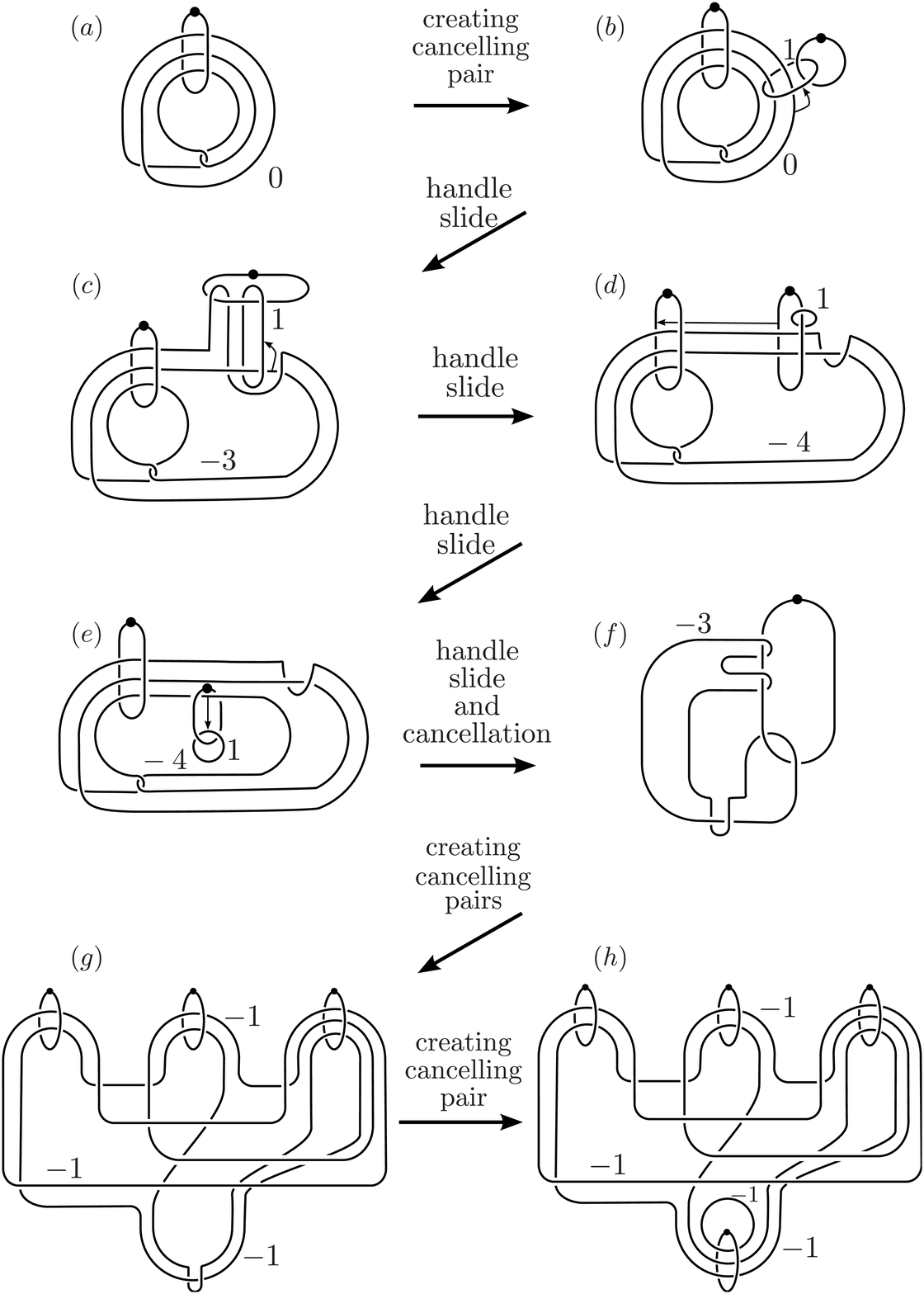}
        \caption{}
        \label{Proof_Akbulut_cork_PALF}
    \end{figure}

    We now show that $X$ is diffeomorphic to $W_1$. 
    The obvious Kirby diagram for $W_1$ is given by
    Figure \ref{AkbulutCork}.
    We draw it as in Figure \ref{Proof_Akbulut_cork_PALF}$(a)$,
    and create the cancelling pair to get Figure
    \ref{Proof_Akbulut_cork_PALF}$(b)$.
    We slide the $0$-framed $2$-handle over a $1$-framed $2$-handle
    to get Figure \ref{Proof_Akbulut_cork_PALF}$(c)$.
    We get Figure \ref{Proof_Akbulut_cork_PALF}$(d)$
    by sliding the $-3$-framed $2$-handle over a $1$-framed $2$-handle.
    By $1$-handle slide, 
    we get Figure \ref{Proof_Akbulut_cork_PALF}$(e)$.
    We slide the $-4$-framed $2$-handle over a $1$-framed $2$-handle,
    and erase a cancelling $1$-handle/$2$-handle pairs
    to get Figure \ref{Proof_Akbulut_cork_PALF}$(f)$.
    We create the cancelling pairs 
    to get Figure \ref{Proof_Akbulut_cork_PALF}$(g)$.
    We get Figure \ref{Proof_Akbulut_cork_PALF}$(h)$ by
    by creating the cancelling pairs.

    The Kirby diagram for $X$ corresponding to the monodromy
    representation $t_{\alpha_4}, \ldots, t_{\alpha_1}$ is
    given by Figure \ref{Proof_Akbulut_cork_PALF}$(h)$.

    Therefore we conclude that $X$ is diffeomorphic to $W_1$, 
    which implies the theorem.

    \hfill $\Box$

    \vspace*{5mm}

    {\itshape Proof of Theorem \ref{main2}.}
    Let $F_{C_i(m, 1, 3, 0)}$ $(i = 1, 2)$ be the compact oriented surface of
    genus zero with $-m+5$ boundary components and
    $\beta_1, \ldots, \beta_{-m+5}$ and $\gamma_1, \ldots, \gamma_{-m+5}$
    be the curves on $F_{C_i(m, 1, 3, 0)}$ shown in
    Figure \ref{PALF_of_C1} and Figure \ref{PALF_of_C2}, respectively.
    Let $g_1:X_{C_1(m, 1, 3, 0)} \rightarrow D^2$
    (resp. $g_2:X_{C_2(m, 1, 3, 0)} \rightarrow D^2$) be
    a Lefschetz fibration over $D^2$ 
    with monodromy representation $(t_{\beta_{-m+5}},\ldots ,t_{\beta_1})$
    (resp. $(t_{\gamma_{-m+5}}, \ldots, t_{\gamma_1})$). 
    Since each curve $\beta_i$ (resp. $\gamma_i$) is
    homologically non-trivial on $F_{C_1(m, 1, 3, 0)}$ (resp. $F_{C_2(m, 1, 3, 0)}$), 
    we see that $g_1$ (resp. $g_2$) is a PALF with fiber
    $F_{C_1(m, 1, 3, 0)}$ (resp. $F_{C_2(m, 1, 3, 0)}$).

\begin{figure}[htbp]
        \centering
        \includegraphics[width=120mm]{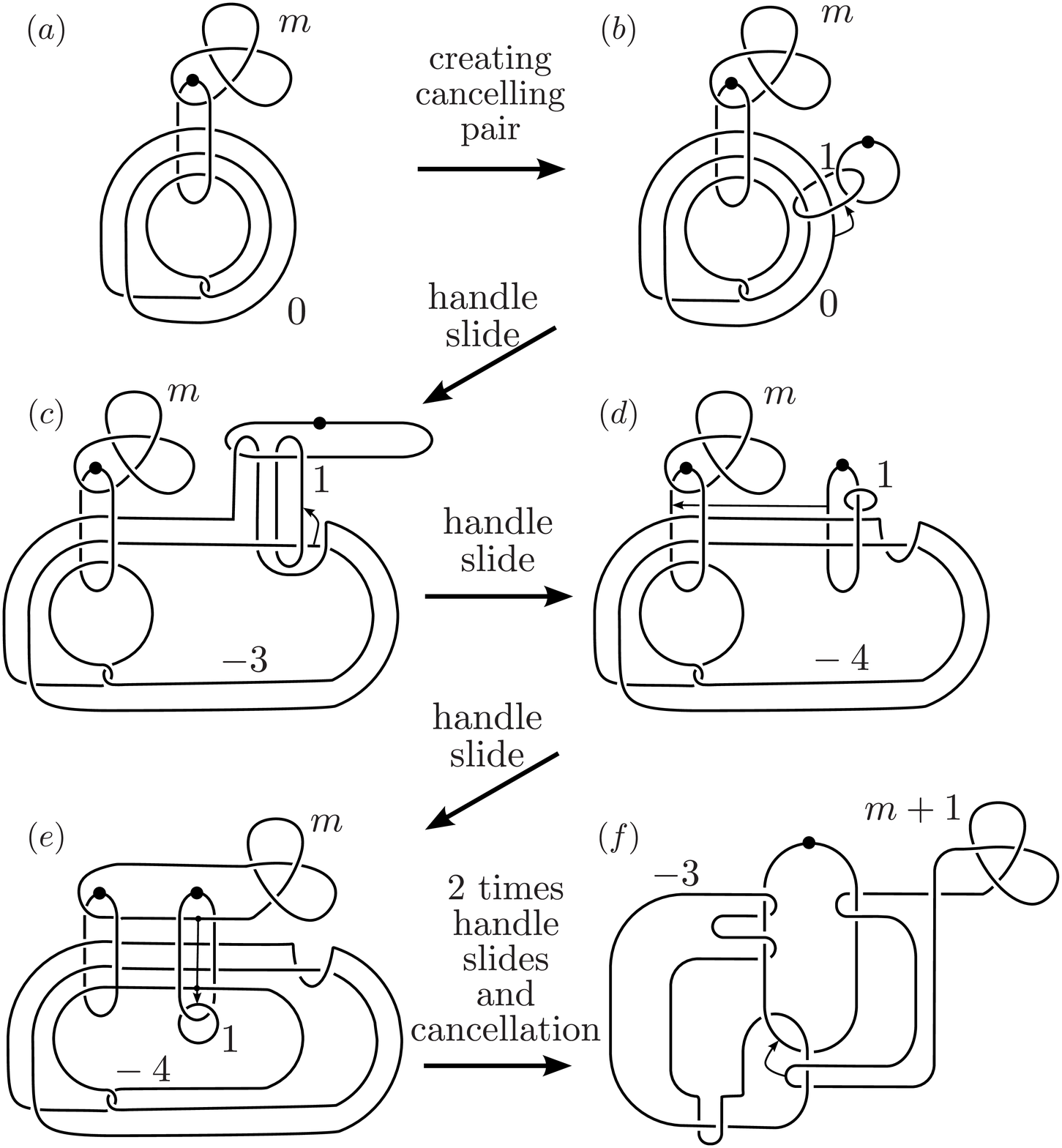}
        \caption{}
        \label{Proof_C1_PALF_1}
    \end{figure}

\begin{figure}[htbp]
        \centering
        \includegraphics[width=140mm]{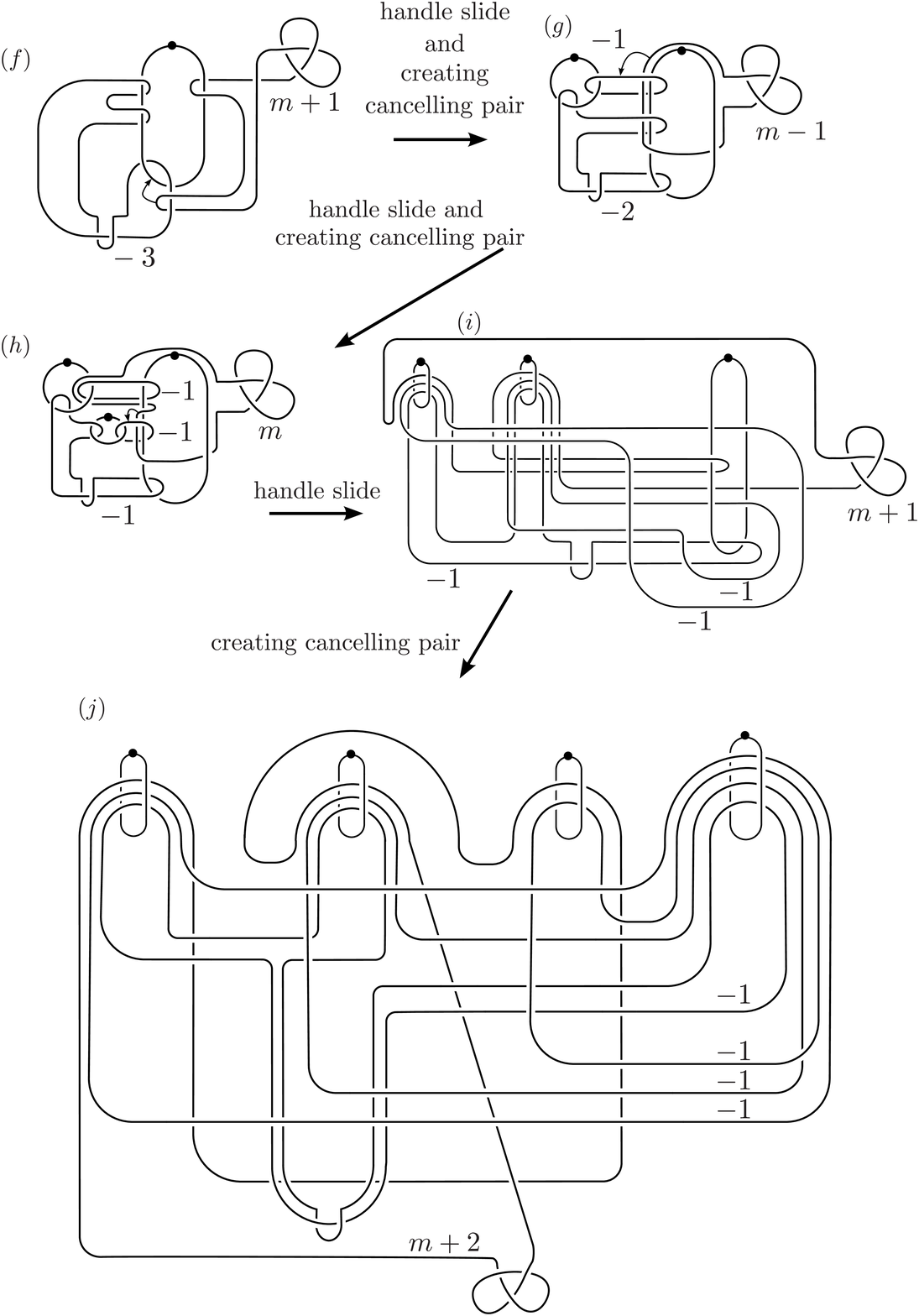}
        \caption{}
        \label{Proof_C1_PALF_2}
    \end{figure}

\begin{figure}[htbp]
        \centering
        \includegraphics[width=120mm]{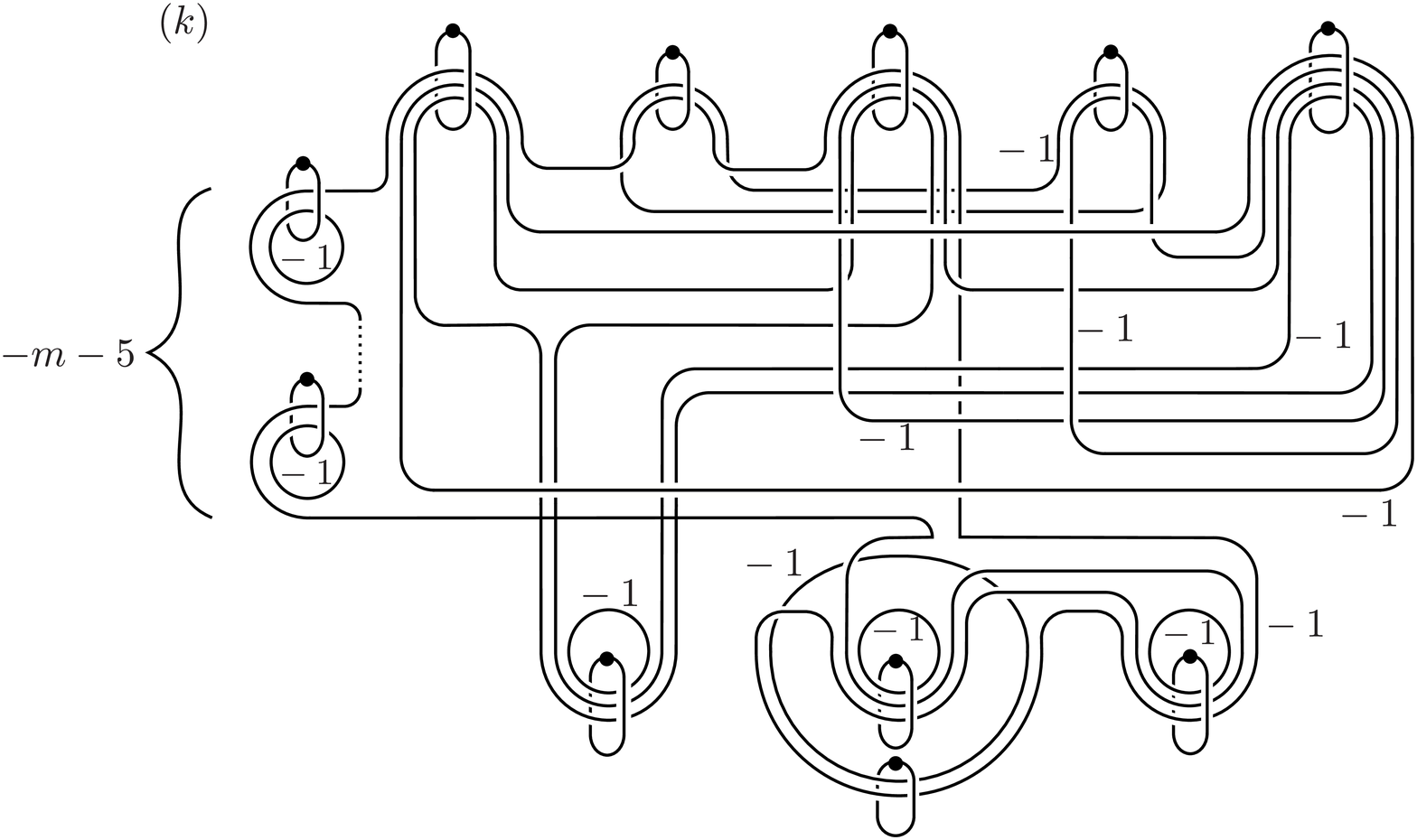}
        \caption{}
        \label{Proof_C1_PALF_3}
    \end{figure}

    The Kirby diagram for $C_1(m, 1, 3, 0)$ is given by
    Figure \ref{Kirby_diag_of_C1}.
    We get Figure \ref{Proof_C1_PALF_1}$(a)$ by isotopy.
    We create a cancelling pair to get Figure \ref{Proof_C1_PALF_1}$(b)$.
    We slide the $0$-framed $2$-handle over a $1$-framed $2$-handle to get
    Figure \ref{Proof_C1_PALF_1}$(c)$.
    By sliding the $-3$-framed $2$-handle over a $1$-framed $2$-handle,
    we get Figure \ref{Proof_C1_PALF_1}$(d)$.
    We get Figure \ref{Proof_C1_PALF_1}$(e)$ by $1$-handle slide.
    We obtain Figure \ref{Proof_C1_PALF_1}$(f)$ by handle slides and
    cancellation.
    By handle slide and creating cancelling pair,
    we get Figure \ref{Proof_C1_PALF_2}$(g)$.
    We get Figure \ref{Proof_C1_PALF_2}$(h)$ by handle slide and
    creating cancelling pair.
    We slide the $m$-framed $2$-handle over $-1$-framed $2$-handle
    to get Figure \ref{Proof_C1_PALF_2}$(i)$.
    In Figure \ref{Proof_C1_PALF_2}$(i)$, creating cancelling pair
    gives Figure \ref{Proof_C1_PALF_2}$(j)$.
    We create cancelling pairs to get Figure \ref{Proof_C1_PALF_3}$(k)$.

    The Kirby diagram for $X_{C_1(m, 1, 3, 0)}$ corresponding to the
    monodromy representation $t_{\beta_{-m+5}}, \ldots, t_{\beta_1}$ is
    given by \ref{Proof_C1_PALF_3}$(k)$.

\begin{figure}[htbp]
        \centering
        \includegraphics[width=100mm]{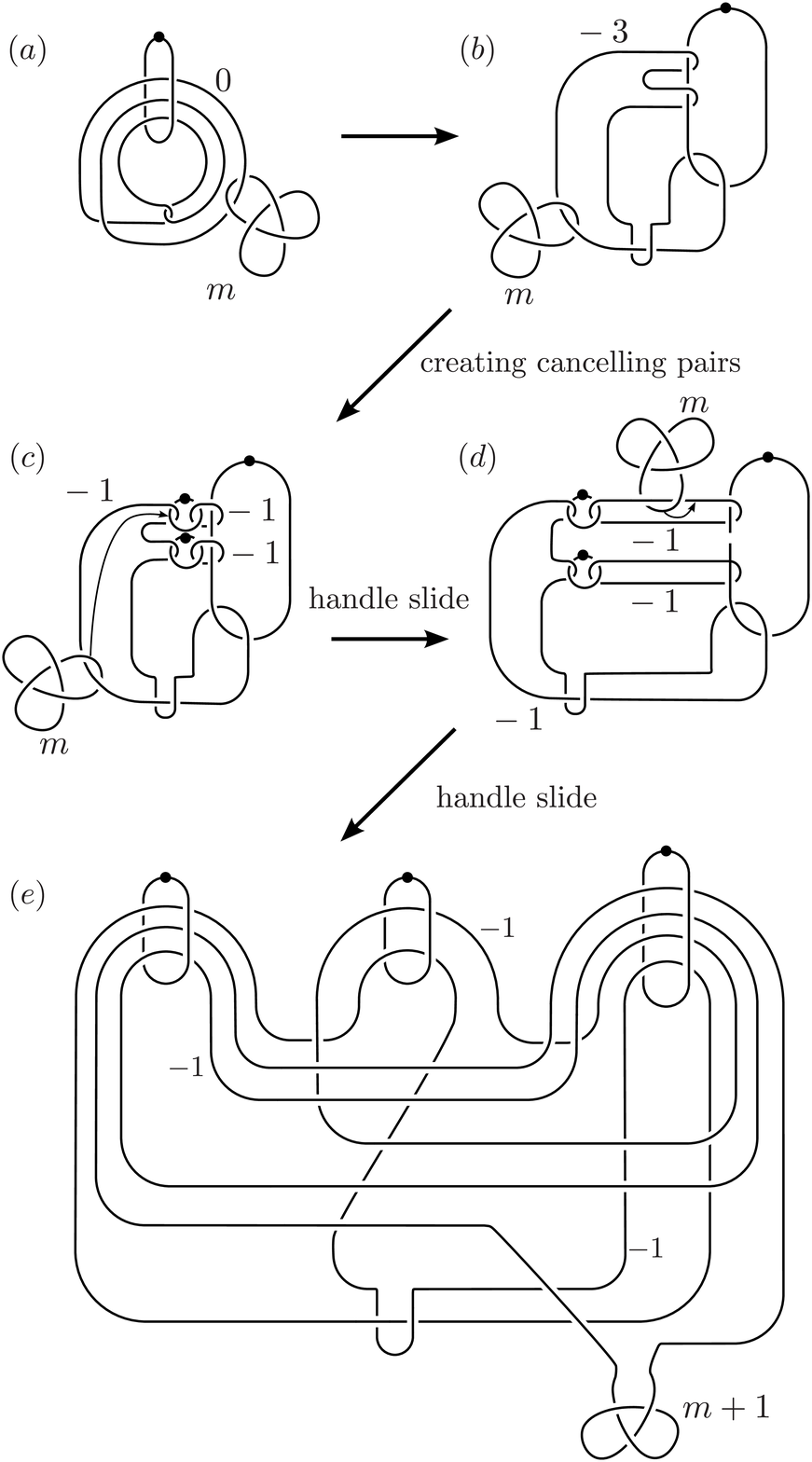}
        \caption{}
        \label{Proof_C2_PALF_1}
    \end{figure}

\begin{figure}[htbp]
        \centering
        \includegraphics[width=120mm]{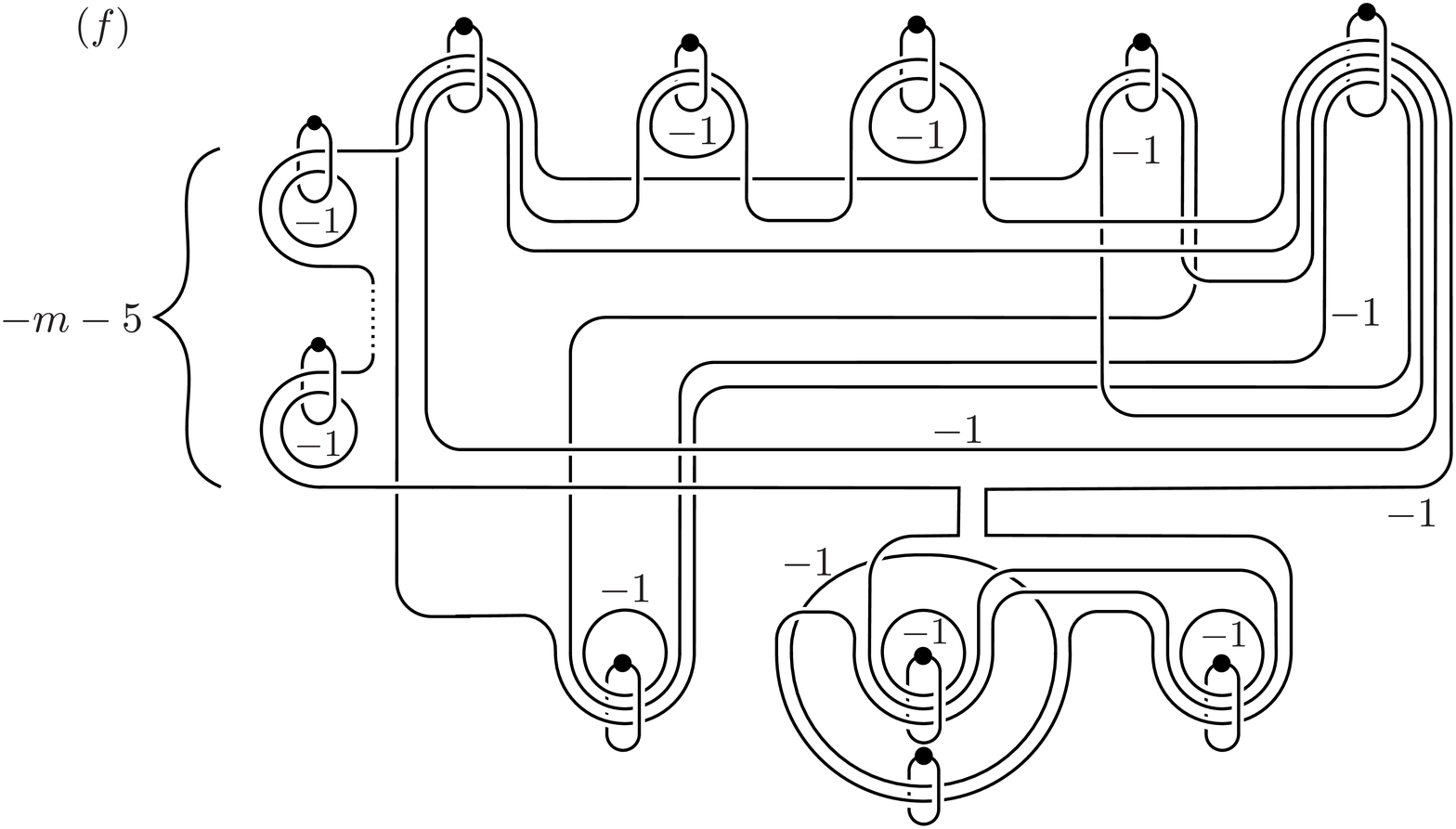}
        \caption{}
        \label{Proof_C2_PALF_2}
    \end{figure}

    The Kirby diagram for $C_2(m, 1, 3, 0)$ is given by
    Figure \ref{Kirby_diag_of_C2}.
    We get Figure \ref{Proof_C2_PALF_1}$(a)$ by isotopy.
    We obtain Figure \ref{Proof_C2_PALF_1}$(b)$ by handle slides and
    cancellation (see the proof of Theorem \ref{main1}).
    By creating cancelling pairs, we get Figure \ref{Proof_C2_PALF_1}$(c)$.
    We slide the $m$-framed $2$-handle under a $1$-handle to get
    Figure \ref{Proof_C2_PALF_1}$(d)$.
    In Figure \ref{Proof_C2_PALF_1}$(d)$, handle slides gives
    Figure \ref{Proof_C2_PALF_1}$(e)$.
    We create cancelling pairs to get Figure \ref{Proof_C2_PALF_2}$(f)$.

    The Kirby diagram for $X_{C_2(m, 1, 3, 0)}$ corresponding to
    the monodromy representation $t_{\gamma_{-m+5}}, \ldots, t_{\gamma_1}$
    is given by Figure \ref{Proof_C2_PALF_2}$(f)$.

    Therefore we conclude that
    the manifolds $C_1(m, 1, 3, 0)$ and $C_2(m, 1, 3, 0)$ admit
    genus zero PALFs.
    \hfill $\Box$

    %%%%%%%%%%%%%%%%%%%%%%%%%%%%%%%%%%%%%%%%%%%%%%%%%%

    %%%%%%%%%%%%%%%%%%%%%%%%%%%%%%%%%%%%%%%%%%%%%%%%%%

    \clearpage

    \renewcommand{\abstractname}{Acknowledgements}
    \begin{abstract}
        The author would like to thank his adviser
        Hisaaki Endo for his helpful comments and
        his encouragement.
        The author wishes to thank Kouichi Yasui for
        his useful comments.
    \end{abstract}

\end{document}